# TAUBERIAN THEOREMS AND LARGE DEVIATIONS
# N. H. BINGHAM

*Abstract.* The link between Tauberian theorems and large deviations is surveyed, with particular reference to regular variation.

*1. Introduction.*

To this author, limit theorems in probability theory have always been irresistibly attractive, as well as centrally important. As Kolmogorov put it in his introduction to Bernoulli's book, "The cognitive value of probability theory lies in the establishment of strict regularities resulting from the combined effects of mass random phenomena. The very notion of mathematical probability would have been fruitless if it were not realized as the frequency of a certain result under repeated experimentation. That is why the works of Pascal and Fermat can be viewed as only the prehistory of probability, while its true history begins with J. Bernoulli's law of large numbers" ([22]; [37], 874-5). The weak law of large numbers gives us a limiting frequency in the presence of a mean; its refinement, the central limit theorem, gives us a limiting distribution in the further presence of a variance. Both phenomena can occur more generally. The necessary and sufficient condition is known in each case, and in each case involves regular variation (for details see e.g. [4]). Thus probabilists are led inescapably towards regular variation (for which see e.g. [5], particularly Chapter 8). Furthermore, the role of the characteristic function, or Fourier transform, in translating averaging, or addition of independent random variables, into multiplication rather than convolution leads one to work with transforms rather than with distributions. But one needs eventually to translate the information thus obtained from transform to distribution language. This is the role of Tauberian theory. For a magisterial treatment of the whole subject, see [24]; for regular variation, see Chapter IV there, or [5] Chapter 4.

For background, we mention here what Tauberian theorems are. When handling limits, or asymptotic properties, we often need to pass from an object to some smoothed version of it. This is usually easy, and is an *Abelian* theorem (prototypes: Abel's continuity theorem for power series, or the implication from $s_n \to c$ to $\sigma_n \to c$, with $\sigma_n$ the arithmetic mean $(s_1 + \ldots + s_n)/n$).

The converse is false in general, but may become true under some additional condition. Such a result is a *Tauberian theorem*, and the additional condition (the 'signature' of the Tauberian theorem) is a *Tauberian condition*. Sometimes one may have a hypothesis on both object and smoothed version, and a conclusion on the object, with *no* additional, or Tauberian, theorem. Such a result is a *Mercerian* theorem (prototype: Mercer's theorem, that $(s_n + \sigma_n)/2 \to c$ implies $s_n \to c$).

Our starting point here is §4.12 of [5], on Tauberian theorems of exponential type. There the treatment closely followed that of Kasahara [19] (leading to three Tauberian theorems, named after de Bruijn, Kohlbecker and Kasahara). At that time, neither Yuji Kasahara nor I realized the link between this area and large deviations. This link has been subsequently elaborated, by Kasahara and Kosugi; see [25], [26], [20], [21]. It involves convex analysis and the Gibbs variational principle; for background, see [40], [17], [11], [10]. The theme is exemplified by the relation displayed on the cover of [10]:

$$-\frac{1}{n} \log \int_x e^{-nh} d\theta^n \to \inf_{x \in X} \{h(x) + I(x)\}$$

($I$ is the *rate function*).

For later use, we record here our notation for the *Fenchel dual* (or Legendre-Fenchel dual, or conjugate) of a convex function $f$:

$$f^*(x) := sup_y\{xy - f(y)\}$$

([36] Sections 7, 12). Then $f^*$ is convex, proper ($< +\infty$ somewhere, $> -\infty$ everywhere) iff $f$ is proper; the closure (greatest lower semi-continuous minorant if $f$ is nowhere $-\infty$, identically $-\infty$ otherwise) cl $f$ has the same dual as $f$, $f^*$ is closed (equal to its closure), and $f^{**} = $ cl $f$.

We shall only be concerned here with the case when all moments $\mu_n := E[X^n]$ of $X$ are finite (in this case the moment-generating function (mgf) $M(t) := E[e^{tX}] = \sum_{n=0}^{\infty} \mu_n t^n/n!$ is defined, though its radius of convergence may be zero). For $X$ non-negative (which we may suppose here without loss) the non-integer moments are then finite also, and one has the option of working with the *Mellin transform* (or Mellin-Stieltjes transform) $\mathcal{F}(s) := E[X^s]$. (Of course, we may pass between the mgf and the Mellin transform by an exponential, or logarithmic, change of random variable, but this may not be probabilistically natural.) Just as the cf is the natural tool for handling *sums* of independent random variables, the Mellin transform is the natural tool for handling *products*. For a systematic treatment, see [41]; for some applications, see Bingham and Teugels [6]. We return to Mellin transforms in Section 4.

2. Entire characteristic function

We know from the work of Paley and Wiener [34] that the proper setting for the Fourier transform is the complex plane (or domain, in their terminology). Thus the traditional 'i' in the definition $\phi(t) := \int e^{itx} dF(x)$ of the characteristic function (cf) $\phi$ of a probability distribution $F$ is merely conventional, and may be omitted if we intend to let the argument $t$ take complex values. One then needs to distinguish cases according to the nature of the cf, or mgf, as a function of a complex variable. It may be entire, analytic but not entire or non-analytic, according as the radius of convergence is infinite, positive and finite, or zero. We turn now to the entire case, for which see e.g. [29], and [5] Section 7.6 (Levin-Pfluger theory: functions of completely regular growth). For the analytic case, see [29] Ch. 7, [30] Ch. 7. For the non-analytic case, see Section 3 below.

An entire cf (of a positive random variable $X$, with distribution function $F$) has order $\rho \in [1, \infty)$ ([29], Th. 7.1.3). If for large $t$ $\log E(e^{tX})/t^\rho$ has liminf $a$ and limsup $b$ with $0 < a \leq b < \infty$, we say that $X$ has *very regular growth* ([38] II.7). Davies [9] shows that in this case one can find upper and lower bounds on

$$-\log P(X > x)/x^{\rho/(\rho-1)},$$

in terms of $a$ and $b$. If the limit exists, and $a = b = \alpha$, say, then we say that $X$ has *completely regular growth*. The Levin-Pfluger theory of [28], [5] Section 7.6 then applies, and $-\log P(X > x)/x^{\rho/(\rho-1)}$ is regularly varying with index $p = \rho/(\rho - 1)$.

This last result is extended to the case of general regular variation by Geluk [14]. Writing $R_p$ for the class of functions regularly varying with index $p$, he links

$$\liminf_{x \to \infty} -\log P(X > x)/\phi(x) = 1$$

for $\phi \in R_p$ with

$$\limsup_{x \to \infty} [E(X^n)]^{1/n} \psi(n)/n = q^{1/q} e^{-1/p},$$

where $p > 1$, $q > 1$ are conjugate (in the usual sense of $L_p$ spaces) and $\psi$ is increasing and in $R_{1/q}$. Then

$$\psi^\rightarrow(x) \sim \phi^*(x),$$

where one has the inverse function of $\psi$ on the left and the Fenchel dual of $\phi$ on the right, so both sides are in $R_q$ (see e.g. [5] Section 1.8.4).

*Example: The Mittag-Leffler distribution.*

For $\alpha \in (0,1)$, the *Mittag-Leffler law* with *index* (or parameter) $\alpha$ has mgf

$$ML_\alpha(t) := \sum_{n=0}^{\infty} t^n / \Gamma(1+n\alpha);$$

for uses in probability theory, see e.g. [3], [12] XIII.8, [5] Section 8.11. Thus

$$\mu_n = n!/\Gamma(1+n\alpha) = \Gamma(n)/(\alpha\Gamma(n\alpha)),$$

and by Stirling's formula one finds, in the notation above,

$$\psi(x) = x^\alpha, \qquad \psi^{\rightarrow}(x) \sim \phi^*(x) \sim x^{1/\alpha},$$

and so

$$-\log P(X > x) \sim \phi(x) \sim (1-\alpha)\alpha^{\alpha/(1-\alpha)} x^{1/(1-\alpha)}.$$

This result is derived twice in [6], p. 350, 352, and again in [5], Th. 8.1.12, as an instance of Kasahara's Tauberian theorem. For the mgf itself, one has

$$ML_\alpha(t) \sim \exp\{t^{1/\alpha}\}/\alpha$$

(Cartwright [8], p. 50). The Mittag-Leffler law of index $\alpha$ is closely connected to the stable subordinator (increasing stable process) of index $\alpha$, and the result above translates into a tail estimate at zero for $X = X(1)$ in the stable subordinator:

$$-\log P(X < x) \sim (1-\alpha)\alpha^{\alpha/(1-\alpha)}/x^{\alpha/(1-\alpha)} \qquad (x \downarrow 0);$$

see [6] 349-350, [5] Section 8.3.4.

The limiting case $p = 1$ relates to entire cfs of infinite order. This is a boundary case of the above, and may be handled by using second-order theory of regular variation (de Haan theory: [5] Chapter 3). Let $\phi$ be a positive function with $\phi(t)/t \to \infty$ as $t \to \infty$, and similarly for $a$. Suppose the second-order condition

$$\frac{\phi(tx) - x\phi(t)}{ta(t)} \to x \log x \qquad (t \to \infty) \qquad \forall x > 0$$

holds. Then if the tail function is everywhere positive and the cf is entire, the following are equivalent:

$$\liminf_{x \to \infty} \frac{-\log P(|X| > x) - \phi(x)}{xa(x)} = 0$$

for some $\phi$, $a$ as above,

$$\liminf_{n \to \infty} \frac{[E(|X|^n)]^{-1/n} n!^{1/n} - \psi(n)}{\alpha(n)} = 0$$

for some $\psi$ and some $\alpha \in \Pi$ (the de Haan class $\Pi$; see e.g. [5] Ch. 3). Then

$$\phi^*(x) \sim (ex)^{-1} \int_0^x \psi^\rightarrow(u) du$$

and $\alpha \sim a$. We refer to Geluk [15], Geluk, de Haan and Stadtmüller [16] for details.

3. *Non-analytic characteristic function*

A cf is defined in the first instance on a line in the complex plane (the real axis if one includes the '$i$' in the definition, the imaginary axis if, as here, one does not), but can be continued analytically to a strip in the complex plane (which may be the whole plane, a half-plane, a strip of positive finite width or the imaginary axis). A cf is analytic iff all moments $\mu_n$ are finite, and $\limsup_{n \to \infty}[|\mu_n|/n!]^{1/n} = 1/R$ is finite ($R$ is then the radius of convergence of the Taylor series of the cf, or mgf). Equivalently, it is analytic iff for some $R > 0$ the tail function $T(x) := P(|X| > x)$ is $O(e^{-rx})$ for all $r < R$, and then the strip of regularity contains the strip $|\operatorname{Im} z| < R$ ([29], Th. 7.2.1). We proceed to consider results analogous to those above, but with the tails heavy enough to cause non-analyticity – that is, with the log-tail decaying more slowly than a multiple of $x$.

A recent result of this kind, which may serve as a prototype, is due to König and Mörters ([23], Lemma 2.3): for $X$ non-negative, $p > 1$ and $c$ real,

$$\lim_{n \to \infty} \frac{1}{n} \log E[X^n/(n!)^p] = -c$$

implies

$$\lim_{x \to \infty} x^{-1/p} \log P(X > x) = -pe^{c/p},$$

with a similar implication with limsup in place of limit. (In their setting, $p = 1, 2, \ldots$, and they consider the intersection of $p$ independent Brownian motions in $d$ dimensions. We exclude here the case $p = 1$, which falls under analyticity above, and is not relevant to their intersection problem, but there is no need here to restrict to $p$ integer.) Their proof is simple, using only Markov's inequality, Stirling's formula and reduction of the limit to a suitable sequence. On the other hand, the setting of their main results makes heavy use of large-deviation theory.

4. *Complements*

1. *Laplace's method.* The essence of Tauberian theorems of exponential type, and of large deviations, is that growth is dominated by behaviour in the neighbourhood of the maximum. This is simply a reflection of how overwhelmingly rapid exponential growth (or decay) is. The prototype in analysis is Laplace's method (for which see e.g. de Bruijn [7]), which goes back to Laplace's *Théorie analytique des probabilités* of 1820 and of which the most telling instance is Stirling's formula of 1730. Its applicability in statistical mechanics rests ultimately on the enormous magnitude of Avogadro's constant, c. $6 \times 10^{23}$.

2. *The Wiman-Valiron method.* In the discrete analogue of Laplace's method (as in [6], 346), one examines the growth of an entire function (an entire cf in our case) by localizing at the term of maximum modulus in its Taylor series. For surveys of the Wiman-Valiron method, see Hayman [18], Fenton [13].

3. *Esscher transforms.* For $F$ a law with mgf $M(t)$, the law $F_t(dx) := e^{tx} dF(x)/M(t)$ is called the *associated* law by Feller [12], XVI.7, or the *exponentially tilted* law. Its use goes back to F. Esscher in 1932, and later to Cramér in 1938; it is accordingly called the *Esscher transform*. It is widely used in ruin theory and insurance and actuarial mathematics; see e.g. Mikosch [31], Section 4.2. In brief, one may use it to change at will the point in whose neighbourhood one wishes to probe the asymptotic behaviour.

4. *Mellin transforms.* As mentioned earlier, with all moments finite as here, one may consider the Mellin transform $\mathcal{F}(s) := E[X^s]$ of $X$, and seek to link its behaviour for large $s$ with the tail behaviour of $X$. One early example is the work of Wagner [39], also considered in [6] (the beginning of my collaboration with Jef Teugels, which led to [5], the origin of Fenchel duality in regular variation, and one source, together with [1], of smooth variation – see [5] Section 1.8).

5. *Branching processes.* Tauberian theorems of exponential type are also useful in branching processes; see e.g. [2]. Again, the theme is that when populations grow, they grow so rapidly that the contribution near the maximum is all that counts.

6. *Random potentials.* We close by mentioning the wealth of related literature on random Schödinger operators in models of disordered media, etc. For a recent survey, with many references, see Leschke et al. [27]. It is no accident that Tauberian theorems of exponential type were also discovered in the theoretical physics literature, as the Minlos-Povsner Tauberian theorem [32], [33]. As Cindy says, all good mathematicians are also interested in physics.

*Postscript.* It is a pleasure to record a personal debt to Cindy Greenwood here. The long Chapter 4 in [5] (64 pages), entitled *Abelian and Tauberian theorems*, originally covered Mercerian theorems also. Only after repeated and urgent pleas from Cindy did I detach the Mercerian component, which became Chapter 5 (25 pages – with Chapter V in Pitt [35], incidentally, the only textbook chapters I know, though see also Korevaar [24] IV.13, 14, Paley and Wiener [34] IV.18). She pointed out with remorseless clarity that separated, the Mercerian material would be impossible to miss, but unseparated risked sinking without trace. I wish I had had the courage of her convictions and split the Abelian and Tauberian parts also.

*Acknowledgement.* I thank Peter Mörters for rekindling my interest in this subject.

Mathematics Department, Imperial College, London SW7 2AZ
n.bingham@ic.ac.uk  nick.bingham@btinternet.com